\documentclass[11pt, a4paper, draft]{amsart}

\usepackage{amsfonts}
\usepackage{amsmath}
\usepackage{amsthm}
\usepackage{amscd}
\usepackage{euscript}

\newcommand{\ZZ}{{\mathbb Z}}

\newcommand{\CC}{{\mathbb C}}

\newcommand{\dual}{{*}}

\newcommand{\pp}{{\mathbf P}}
\newcommand{\PP}{{\mathbb P}}
\newcommand{\OO}{{\EuScript O}}

\newcommand{\PPn}{{\PP^n}}
\newcommand{\PPm}{{\PP^m}}
\newcommand{\PPt}{{{\widetilde\PP}^m}}
\newcommand{\PPV}{{\PP(V)}}
\newcommand{\OOn}{{\OO_\PPn}}

\newcommand{\netp}{{NE(\PPt)}}
\newcommand{\calm}{{\mathfrak M}}

\newcommand{\hdue}{{\HH_2(\PPV)}}

\DeclareMathOperator{\PSL}{\PP SL}
\DeclareMathOperator{\HH}{H}

\DeclareMathOperator{\Ext}{Ext}
\DeclareMathOperator{\im}{Im}

\DeclareMathOperator{\virtdim}{{virtdim}}

\newcommand{\iso}{\simeq}

\newcommand{\qp}{*}

\newtheorem{thm}{Theorem}[section]
\newtheorem{lemma}[thm]{Lemma}
\newtheorem{cor}[thm]{Corollary}
\newtheorem{prop}[thm]{Proposition}

\theoremstyle{definition}

\theoremstyle{remark}

\begin{document}
\title[On the quantum cohomology of blow-ups]{On the quantum
  cohomology of blow-ups \\ of projective spaces along linear
  subspaces.}

\author{Marco Maggesi}
\email{maggesi@math.unifi.it}
\urladdr{http://www.math.unifi.it/\~{}maggesi/}
\begin{abstract}
  We give an explicit presentation with generators and
  relations of the quantum cohomology ring of the blow-up of a
  projective space along a linear subspace.
\end{abstract}

\maketitle

\section{Introduction.}

Let $\PPm$ be the complex projective space, $\Lambda\subset\PPm$ a
linear subspace of dimension $p$ and $\alpha \colon \PPt \to \PPm$ the
blow-up of $\PPm$ along $\Lambda$.  Let $k$ be the hyperplane class on
$\PPm$ and $\eta$ the exceptional divisor on $\PPt$.  The aim of the
present paper is to show a way to compute the quantum cohomology ring
$\HH_Q(\PPt)$ of $\PPt$.  The (classical) cohomology ring of $\PPt$
can be expressed as (compare with lemma \ref{t:c-ring})

\begin{equation}
  \HH^*(\PPt) = \ZZ[k,\eta] / (g_1, g_2),
\end{equation}
where the two relations are
\begin{equation}
  \label{e:c-rel-g}
  g_1(k,\eta)=(k-\eta)^{m-p},\qquad g_2(k,\eta)=k^{p+1}\eta.
\end{equation}

The main result of this work is the following
\begin{thm}[Main theorem]
  \label{t:main}
  Suppose that $2p+3<m$.  Then the quantum cohomology ring of $\PPt$
  can be expressed as
  \begin{equation}
    \label{e:q-ring}
    \HH^*_Q(\PPt) = \ZZ[k,\eta] / (\tilde g_1, \tilde g_2),
  \end{equation}
  where the two relations are
  \begin{equation}
    \label{e:q-rel-g}
    \tilde g_1(k,\eta)=(k-\eta)^{m-p}-\eta,
    \qquad \tilde g_2(k,\eta)=k^{p+1}\eta-1.
  \end{equation}
\end{thm}

The blow-up $\PPt$ can be regarded as a projective bundle on a
projective space (proposition \ref{t:isom}).  The quantum cohomology
of projective bundles on projective spaces was studied by Qin and Ruan
in \cite{qr}.  The relevant material from their work will be enclosed
here for completeness.  Other studies on quantum cohomology of
blow-ups of projective spaces can be found in \cite{bg}, \cite{ga1}
and \cite{ga2}.

I wish to thank professor V.~Ancona for many useful discussions and
for his encouragement.

\section{$\PPt$ as projective bundle.}
Let $n:=m-p-1$, $ r:=p+2$ and $V$ be the rank-$r$ vector bundle on
$\PPn$ given by
\begin{equation}
  \label{d:V}
  V:=\OOn(1)^{\oplus r-1}\oplus\OOn(2).
\end{equation}
and consider the associate $\PP^{r-1}$-bundle $\pi \colon \PPV \to \PPn$.

\begin{prop}
  \label{t:isom}
  The two varieties $\PPt$ and $\PPV$ are isomorphic.
\end{prop}

\begin{proof}
  Consider $\PPm$ as the projective space $\pp(U)$ of lines of a
  vector space $U$ of dimension $m+1$ and $\Lambda=\pp(U_0)$, where
  $U_0$ is a subspace of $U$ of dimension $p+1$ (We will use ``$\PP$''
  for projective spaces of hyperplanes and ``$\pp$'' for projective
  spaces of lines).  Two subspaces $L$, $W$ of $U$ of dimension $1$
  and $p+2$ determine points \( [L] \in \pp(U) \) and \( [W/U_0] \in
  \pp(U/U_0) \) respectively.
  The blow-up $\PPt$ can be identified with the closed incidence
  subvariety \( Y \subset \pp(U) \times \pp(U/U_0) \) defined by \( Y
  := \big\{ ([L],[W/U_0]) \colon L \subset W \big\} \) and \( \alpha
  \colon \PPt \to \PPm \) is the projection on the first factor
  restricted to $Y$:
  \begin{equation}
    \label{dg:Y}
    \begin{CD}
      Y @>>> \pp(U/U_0)\\
      @VVV @.\\
      \pp(U) @.  {}
    \end{CD}
    \qquad\qquad
    \begin{CD}
      ([L],[W/U_0]) @>>> [W/U_0] \\
      @VVV    @.\\
      [L] @.  {}
    \end{CD}
  \end{equation}
  Note that \( Y \to \pp(U/U_0) \) is a projective bundle \( \pp(E)
  \to \pp(U/U_0) \), where $E$ is the vector bundle \( E = \{(w,[W])
  \in U \times \pp(U/U_0)\} \).  We want to show that $\pp(E)$ and
  $\PPV$ are isomorphic.
  
  The exact sequence $0 \to U_0 \to W \to U/U_0 \to 0$ gives rise to
  an exact sequence of vector bundles
  \begin{equation}
    \label{s:E}
    0 \to U_0\otimes\OOn \to E \to \OOn(-1) \to 0.
  \end{equation}
  The above sequence splits (\( \Ext^1 (\OOn(-1), U_0\otimes\OOn) = 0
  \)), so \( E = (U_0\otimes\OOn) \oplus \OOn(-1) \).  Since \( V =
  E^\dual(1) \), then
  \begin{math}
    \PP(V) \iso \PP(E^\dual(1)) \iso \PP(E^\dual) \iso \pp(E).
  \end{math}
\end{proof}

The exceptional locus \( \alpha^{-1}(\Lambda) \simeq \Lambda \times
\PPn \) of the blow-up \( \alpha \colon \PPt \to \PPm \) corresponds
to the trivial sub-bundle \( \PP(U_0 \otimes \OO_\PPn(1)) \) of \( \pi
\colon \PPV \to \PPn \) under the isomorphism of proposition
\ref{t:isom}.  The two restrictions \( \alpha|_{\Lambda \times \PPn}
\colon \Lambda \times \PPn \to \Lambda \) and \( \pi|_{\PP(U_0 \otimes
  \OO_\PPn(1))} \colon \PP(U_0 \otimes \OO_\PPn(1)) \to \PPn \) are
the canonical projections of the product \( \PP(U_0) \times \PPn \).

Let $h$ be the hyperplane bundle on $\PPn$ and $\xi=\xi_V$ the
tautological line bundle on $\PPV$.  We will make no distinction
between $h$ and their corresponding pull-back $\pi^*h$ on $\PPV$.
Note that $\xi-h$ and $h$ are nef on $\PPV$.  It is well known that
the classical cohomology ring of $\PPV$ is generated by $h$, $\xi$
with the two relations
\begin{eqnarray}
  \label{e:c-rel-f}
  f_1(h,\xi) &=& h^{n+1}\\
  f_2(h,\xi) &=& (\xi-h)^{r-1}(\xi-2h)
  =\sum_{k=0}^r(-1)^kc_kh^k\xi^{r-k},
\end{eqnarray}
where $c_i=c_i(V)$ are the Chern classes of $V$.

The homology classes \( A_1 := (h^n\xi^{r-2})_* \) and \( A_2 :=
(h^{n-1}\xi^{r-1} - rh^n\xi^{r-2})_* \) form a system of generators
for $\hdue$ ($(-)_*$ stands for the Poincare dual) and
\begin{equation}
  \label{e:duali}
  (\xi-h)(A_1)=1,\quad(\xi-h)(A_2)=0,\quad h(A_1)=0,\quad h(A_2)=1.
\end{equation}

\begin{prop}
  \label{t:fano}
  The variety $\PPV$ is Fano and the homology classes $A_1$, $A_2$ are
  extremal rays for $\netp$.  The canonical divisor is
  $-K=-K_\PPV=r(\xi-h)+nh$.
\end{prop}
\begin{proof}
  In fact, \( -K = r\xi + (n+1-c_1(V))h = r(\xi-h) + nh \).  Let $C$
  be an effective curve and \( [C] = a_1A_1 + a_2A_2 \).  From
  (\ref{e:duali}) we have \( (\xi-h) \cdot [C] = a_1 \) and \( h \cdot
  [C] = a_2 \).  Since $\xi-h$ and $h$ are nef, $a_1\geq 0$ and
  $a_2\geq 0$.  Moreover, if $C$ is not constant, $a_1>0$ or $a_2>0$
  and \( -K \cdot [C] = ra_1 + na_2 > 0 \).  Since $h$ and $\xi-h$ are
  globally generated, $\PPV$ is Fano.
\end{proof}

\begin{lemma}
  \label{t:A1-A2}
  The classes $A_1$ and $A_2$ can be represented by rational connected
  curves as follows:
  \begin{enumerate}
  \item \( A_1 = [\ell] \), where $\ell$ is a line contained in a
    fiber of \( \pi \colon \PPV \to \PPn \).
  \item \( A_2 = [\tilde \ell] \), where \( \tilde \ell := [\{P\}
    \times \ell] \subset \PP(U_0) \times \PPn \subset \PPV \) with $P$
    a point in $\PP(U_0)$ and $\ell$ a line in $\PPn$.
  \end{enumerate}
  Moreover, there are no other means to represent $A_1$ and $A_2$ as
  class of a rational connected curve beside (i) and (ii).
\end{lemma}
\begin{proof}
  \begin{enumerate}
  \item If $\ell$ is a line in a fiber $\pi$, then \( \xi \cdot [\ell]
    = 1 \) and \( h \cdot [\ell] = 0 \), hence $[\ell]=A_1$ by
    (\ref{e:duali}).  Conversely, let $\ell$ be a rational connected
    curve in $\PPV$ such that $A_1=[\ell]$.  The curve $E$ is
    contained in a fiber $\pi^{-1}(Q)$ as \( [\pi(\ell)] = h^n \).
    Since \( \xi \cdot A_2 = 1 \), then \( \xi|_{\pi^{-1}(Q)} \cdot
    [\ell] = 1 \) in $\pi^{-1}(Q)$, that is $\ell$ is a line in
    $\pi^{-1}(Q)$.
    
  \item Now let \( \ell \subset \PPn \) be a line and $P\in\PPn$.
    Note that \( (\xi-h)|_{\PP(U_0\otimes\OO_\PPn(1))} =
    {\xi_{V(-1)}}|_{\PP(U_0\otimes\OO_\PPn)} \) is the pull-back of
    the hyperplane divisor of $\PP(U_0)$ through the canonical
    projection \( \PP(U_0) \times \PPn \to \PP(U_0) \).  Then \(
    (\xi-h) \cdot [\ell] = 0 \) and \( h \cdot [\ell] = 1\) therefore
    $[\ell]=A_2$ by (\ref{e:duali}).  Conversely, let $\tilde \ell$ be
    a rational connected curve in $\PPV$ such that \( A_2 = [\tilde
    \ell] \) and \( \ell = \pi(\tilde \ell) \).  Since \( h \cdot
    [\tilde \ell] = 1 \), then \( \pi|_{\tilde \ell} \colon \tilde
    \ell \to \ell \) is an isomorphism and $\ell$ is a line in $\PPn$.
    The inclusion \( \tilde \ell \hookrightarrow \PPV \) is induced by
    a surjective morphism \( f \colon V|_\ell \to \OO_\ell(t) \),
    where \( t = c_1(\xi|_\ell) = \xi \cdot [\tilde \ell] = 1 \) (see
    \cite{Ha} II Prop.\ 7.12).  Every morphism \( V|_\ell \to
    \OO_\ell(1) \) factors through \( V|_\ell \to U_0 \otimes
    \OO_\ell(1) \to \OO_\ell(1) \); hence $\tilde \ell$ is contained
    in \( \PP(U_0) \times \PPn \) and has the form \( \{P\} \times
    \ell \).
  \end{enumerate}
\end{proof}

So far we have described classes $A_1$ and $A_2$ in the language
of projective bundles, now we are going to describe them in the
language of blow-ups.  Note that the preimage $\alpha^{-1} (\Lambda)$
of $\Lambda$ under the blow-up morphism $\alpha \colon \PPt \to \PPm$
is isomorphic to $\Lambda \times \PPn$ and \(
\alpha|_{\Lambda\times\PPn} \colon \Lambda \times \PPn \to \Lambda \)
is the canonical projection on the first factor.  Let $\tilde \ell_1$
be the strict transform in the blow-up of a line $\ell_1 \subset \PPm$
not included in $\Lambda$ and meeting $\Lambda$ in one point and \(
\tilde \ell_2 := \{P\} \times \ell_2 \subset \Lambda \times \PPn =
\alpha^{-1}(\Lambda) \) where $P$ is a point of $\Lambda$ and $\ell_2$
is a line in $\PPn$.

\begin{lemma}
  Under the isomorphism \( \PPt \simeq \PPV \), homology classes \(
  B_1 = [\tilde \ell_1] \) and \( B_2 = [\tilde \ell_2] \) correspond
  respectively to $A_1$ and $A_2$.
\end{lemma}
\begin{proof}
  It is obvious that $B_2$ corresponds to $A_2$ since both are
  described as the class of a curve of the form \( \{P\} \times \ell
  \subset \PP(U_0) \times \PPn \).
  
  With notations from proposition \ref{t:isom}, $\ell_1$ is determined
  by a plane \( S_{\ell_1} \subset U \) such that \( S_{\ell_1} \cap
  U_0 \) has dimension $1$ and \( \tilde \ell_1 \) can be described as
  the subset of $Y$ of points of the form \( ([L],[S_{\ell_1}]) \) for
  all subspaces \( L \subset S_{\ell_1} \) of dimension $1$.  Then
  $\tilde \ell_1$ can be regarded as a line in the fiber of the
  projective bundle \( \pi \colon Y \to \pp(U/U_0) \) over the point
  $[S_{\ell_1}/U_0]$, that is $B_1$ corresponds to $A_1$.
\end{proof}

\begin{cor}
  \label{t:c-ring}
  Under the isomorphism \( \PPt \iso\PPV \), cohomology classes $k$,
  $\eta$ correspond respectively to $\xi-h$, $\xi-2h$.  The classical
  cohomology ring of $\PPt$ is given by
  \begin{equation}
    \label{e:c-ring}
    \HH^*(\PPt,\ZZ)=\ZZ[k,\eta] / (g_1,g_2),
  \end{equation}
  where $g_1(k,\eta)=(k-\eta)^{n+1}=(k-\eta)^{m-p}$,
  $g_2(k,\eta)=k^{r-1}\eta=k^{p+1}\eta$.
\end{cor}
\begin{proof}
  The statement follows easily from (\ref{e:duali}), from equations
  \begin{equation}
    k \cdot B_1 = 1, \quad
    k \cdot B_2 = 0, \quad
    \eta \cdot B_1 = 1, \quad
    \eta \cdot B_2 = -1
  \end{equation}
  and from the previous lemma.
\end{proof}

\section{Study of $\calm(A,0)$.}

Let $\calm(A,0)$ be the moduli space of morphisms $f\colon \PP^1\to
\PPV$ with $[\im f]=A$ where $A$ is a class in $\hdue$.  Recall that
the virtual dimension of $\calm(A,0)$ is
\begin{equation}
  \label{e:virtdim}
  \virtdim(\calm(A,0)) = -K \cdot A + n + r - 1.
\end{equation}

\begin{lemma}
  \label{t:MA1}
  The moduli space \( \calm(A_1,0) / \PSL(2,\CC) \) is smooth, compact
  and has expected dimension.
\end{lemma}
\begin{proof}
  By lemma \ref{t:A1-A2}, \( \calm(A_1,0) / \PSL(2,\CC) \simeq G(2,r)
  \times \PPn \), where $G(2,r)$ is the grassmannian of lines in
  $\PP^{r-1}$; then it is smooth, compact and has dimension \( \dim
  G(2,r) + \dim\PPn = n+2r-4 = -K \cdot A_1 + \dim\PPV - 3 \).
\end{proof}

\begin{lemma}
  \label{t:MAt}
  The moduli space \( \calm(A_2,0) / \PSL(2,\CC) \) is smooth, compact
  and has expected dimension.
\end{lemma}
\begin{proof}
  By lemma \ref{t:A1-A2}, \( \calm(A_2,0) / \PSL(2,\CC) \) is a
  $\PP(U_0)$-bundle on $G(2,n+1)$ Then it is smooth, compact and has
  dimension
  \begin{math}
    \dim G(2,n+1) + \dim \PP(U_0)
    = 2n + r - 4
    = -K \cdot A_2 + \dim\PPV - 3.
  \end{math}
\end{proof}

\section{Computation of Gromov-Witten invariants.}

If $A$ belongs to $\hdue$ and $\alpha$, $\beta$, $\gamma$ are classes
in $\HH^*(\PPV)$, \( I_A(\alpha, \beta, \gamma) \) denotes, as usual,
the Gromov-Witten invariant.  If one assume the genericity condition
\begin{itemize}
\item[($\star$)] the moduli space \( \calm (A,0)/\PSL(2,\CC) \) is
  smooth, compact, of expected dimension \( -K_\PPV \cdot A + (n+r-1)
  - 3\)
\end{itemize}
the enumerative interpretation of the Gromov-Witten invariants holds,
that is, if $B$, $C$, $D$ are three sub-varieties of $\PPV$ in general
position representing classes $\alpha$, $\beta$, $\gamma$, then \(
I_A(\alpha, \beta, \gamma) \) is the number of rational curves,
counted with suitable multiplicity, that intersect $B$, $C$ and $D$.

We recall that a sufficient condition for the smoothness of the space
\( \calm(A,0) / \PSL(2,\CC) \) is given by the vanishing of $h^1(N_f)$
for every map \( f \in \calm(A,0) \), where $N_f$ is the normal
bundle.  Moreover, if $A$ can be represented only by irreducible and
reduced curves, \( \calm(A,0) / \PSL(2,\CC) \) is compact.

\begin{lemma}
  \label{t:GW1}
  If $q_1 + q_2 < br$, then
  \begin{math}
    I_{bA_1}(h^{p_1}\xi^{q_1},h^{p_2}\xi^{q_2},\alpha)=0.
  \end{math}
\end{lemma}
\begin{proof}
  Let $B$, $C$, $D$ be three varieties in general position dual to
  $h^{p_1}\xi^{q_1}$, $h^{p_1}\xi^{q_1}$, $h^p\xi^q$ respectively.  It
  is enough to show that there is no effective connected curve $E$
  representing $A$ and intersecting $B$, $C$, $D$.  In fact, the
  genericity condition ($\star$) can be relaxed by assuming (compare
  with \cite{qr} Lemma 3.7):
  \begin{quotation}
    $h^1(N_f)=0$ for every $f \in \calm(A,0)$ such that $\im(f)$
    intersects $B$, $C$, $D$ and there is no reducible or non reduced
    effective (connected) curve $E$ such that $[E]=A$ and $E$
    intersects $B$, $C$, $D$.
  \end{quotation}
  We can choose $\alpha$ of the form $h^ph^q$ with $q\leq r-1$ and
  $p+p_1+p_2+q+q_1+q_2=n+r-1+br$.  Notice that $\pi(B)$, $\pi(C)$,
  $\pi(D)$ in $\PPn$ are dual to $h^{p_1}$, $h^{p_2}$, $h^p$.  Since
  $A_1$ is an extremal ray, every irreducible component of $E$ is a
  multiple of $A_1$ and hence contained in a fiber of $\pi$.
  Therefore the whole $E$ is contained in a fiber of $\pi$ as $E$ is
  connected.  We deduce that
  $\pi(B)\cap\pi(C)\cap\pi(D)\neq\emptyset$; thus $p_1+p_2+p\leq n$
  and
  \begin{math}
      q_1+q_2 = n+r-1+br-p_1-p_2-p-q \geq br.
  \end{math}
\end{proof}

\begin{lemma}
  \label{t:GW2}
  We have
  \begin{math}
    I_{A_1}(\xi,\xi^{r-1},h^n\xi^{r-1})=1.
  \end{math}
\end{lemma}
\begin{proof}
  By lemma \ref{t:MA1}, $\calm(A_1,0)$ is smooth, compact, of expected
  dimension.  The class $h^n\xi^{r-1}$ is the dual of a point $q$ in
  $\PPV$.  Let $p:=\pi(q)$.  A parametrized curve in $\calm(A_1,0)$
  meeting $q$ has support on a line $\ell\subset\pi^{-1}(q)$.  Since
  $\xi|_{\pi^{-1}(q)}$ is the cohomology class of a hyperplane in
  $\pi^{-1}(q)$, then $I_{A_1}(\xi,\xi^{r-1},h^n\xi^{r-1})=1$.
\end{proof}

\begin{lemma}
  \label{t:GW3}
  Let $\tilde n$ be an integer and \( 1 \leq \tilde n \leq n \).
  Then
  \begin{align}
    \label{e:IA2-0}  
    I_{A_2}(h^{\tilde n},h^{n+1-\tilde n},h^n\xi^{r-2}) &= 1, \\
    \label{e:IA2-1}
    I_{A_2}(h^{\tilde n},h^{n+1-\tilde n},h^{n-1}\xi^{r-1}) &= r-1.
  \end{align}
\end{lemma}

\begin{proof}
  By lemma \ref{t:MAt}, \( \calm(A_2,0) \) is smooth, compact, of
  expected dimension.  As explained in lemma \ref{t:A1-A2}, $A_2$ is
  represented by curves $\tilde \ell$ in \( \PP(U_0) \times \PPn
  \subset \PPV \) and \( \ell := \pi(\tilde \ell) \subset \PPn \) is a
  line.
  
  First we prove equation (\ref{e:IA2-0}).  Let $B$, $C$, $D$ be three
  sub-varieties of $\PPV$ in general position representing $h^{\tilde
    n}$, $h^{n+1-\tilde n}$, $h^n\xi^{r-2}$; then $\pi(B)$, $\pi(C)$
  are linear subspaces of $\PPn$ of dimension $\tilde n$, $\tilde n -
  1$ and $\pi(D)$ is a point.  Assume that \( A_2 = [\tilde \ell] \)
  and $\tilde \ell$ intersects $B$, $C$, $D$.  Then $\ell$ is uniquely
  determined as the only line meeting $\pi(B)$, $\pi(C)$ and $\pi(D)$.
  Since $D$ is a line in the fiber of \( \pi \colon \PPV \to \PPn \),
  then \( D \cap (\PP(U_0) \times \PPn) \) is a set of a single point
  $P\in\PPV$.  It follows that \( \tilde \ell = \{P\} \times \ell
  \subset \PP(U_0) \times \PPn \); hence (\ref{e:IA2-0}) holds.
    
  Now let us show that
  \begin{equation}
    \label{e:IA2-2}
    I_{A_2}(h^{\tilde n}, h^{n+1-\tilde n},
    h^{n-1}\xi^{r-1} + (1-r) h^n\xi^{r-2}) = 0.
  \end{equation}
  By linearity, after the proof of equation (\ref{e:IA2-0}) this is
  equivalent to equation (\ref{e:IA2-1}).  Let $\ell^\prime$ be a
  general line in $\PPn$ and \( f \colon V|_{\ell^\prime} \to
  \OO_{\ell^\prime}(2) \) be a general surjective morphism of vector
  bundles.  Let $\tilde \ell^\prime$ be the image of the induced
  embedding \( \PP(f) \colon \PP^1 = \PP(\OO_{\ell^\prime}(2))
  \hookrightarrow \PP(V_{\ell^\prime}) \).  The curve $\tilde
  \ell^\prime$ represents \( h^{n-1}\xi^{r-1} + (1-r) h^n\xi^{r-2} \)
  since \( \xi \cdot [\tilde \ell] = 2 \) and \( h \cdot [\tilde \ell]
  = 1 \).  Observe that \( \tilde \ell \cap \tilde \ell^\prime =
  \emptyset \), indeed, the two trivial sub-bundles \( \PP(\OOn(2)) \)
  and \( \PP(U_0 \otimes \OOn(1)) \) of $\PPV$, which contain
  $\tilde\ell$ and $\tilde \ell^\prime$ respectively, do not
  intersect.  This implies equation (\ref{e:IA2-2}).
\end{proof}

\section{The quantum cohomology ring of $\PPV$:\\
  proof of the main theorem.}

As we are going to work with two ring structures, classical and
quantum, on the cohomology group $\HH^*(\PPV)$, we have to fix some
notation.  If $\alpha$, $\beta$ are cohomology classes in
$\HH^*(\PPV)$ then \( \alpha^{\qp i} \qp \beta^{\qp j} \) is the
quantum product of $i$ copies of $\alpha$ and $j$ copies of $\beta$.
For the classical product we write $\alpha^i\beta^j$ as usual.  For
any homology class $A\in\hdue$ we denote by \( (\alpha \qp \beta)_A \)
the contribution of $A$ to the product \( \alpha \qp \beta \) defined
as follows.  Let $d$ be the integer \( d := \deg(\alpha) + \deg(\beta)
+ K \cdot A \).  Then \( (\alpha \qp \beta)_A = 0 \) for $d<0$.  For
$d\geq 0$, \( (\alpha \qp \beta)_A \) is the class of degree $d$ in
$\HH^*(\PPV)$ satisfying the property
\begin{equation}
  (\alpha \qp \beta)_A \cdot \gamma_* = I_A(\alpha, \beta, \gamma)
\end{equation}
for every $\gamma\in\HH^*(\PPV)$ with \( \deg\gamma = \dim\PPV - d \),
being \( I_A(\alpha, \beta, \gamma) \) the Gromov-Witten invariant.
The quantum product is given by
\begin{equation}
  \begin{split}
    \alpha \qp \beta
    &= \sum_{A\in\hdue}(\alpha \qp \beta)_A\\
    &= \alpha\beta + (\alpha\qp\beta)_{[E_1]}
       + (\alpha\qp\beta)_{[E_2]} + \dots
  \end{split}
\end{equation}
for some non constant effective curves $[E_1]$, $[E_2]$, \dots

\begin{lemma}
  \label{t:q-rel-f}
  If $r<n$, the two relations
  \begin{equation}
    \label{e:q-rel-f}
    h^{\qp(n+1)}=\xi-2h,\qquad
    (\xi-h)^{\qp(r-1)}\qp(\xi-2h)=1
  \end{equation}
  hold in $\HH^*_Q(\PPV)$.
\end{lemma}
\begin{proof}
  Since \( r = \min(-K \cdot A_1, -K \cdot A_2) \), we have \( h^{\qp
    p} \qp \xi^{\qp q} = h^p \xi^q \) whenever $p+q<r$.  Moreover, for
  $p+q<n$, the contributions to the quantum product \( h^{\qp p} \qp
  \xi^{\qp q} \) can only arise from curves of type $tA_1$.  In
  particular, lemma \ref{t:GW1} implies that \( h^{\qp p} \qp \xi^{\qp
    q} = h^p\xi^q \) for $p<n$ and $q<r$.  By lemma \ref{t:GW2} we
  have \( \xi^{\qp r} = 1 \), then
  \begin{equation*}
    (\xi-h)^{\qp(r-1)}\qp(\xi-2h) = (\xi-h)^{r-1}(\xi-2h)+1 = 1.
  \end{equation*}
  Analogously, by lemma \ref{t:GW3} we have \( (h^{\tilde n} \qp
  h^{n+1-\tilde n})_{A_2} = \xi-2h \), then
  \begin{equation*}
    \begin{split}
      h^{\qp (n+1)}
      &= h^{\qp \tilde n} \qp h^{\qp (n+1-\tilde n)}
      = h^{\tilde n} \qp h^{n+1-\tilde n} \\
      &= h^{n+1}+(h^{\tilde n} \qp h^{n+1-\tilde n})_{A_2} = \xi-2h.
    \end{split}
  \end{equation*}
\end{proof}

We are finally ready to prove the main result of this work.
\begin{proof}[Proof of theorem \ref{t:main}]
  If $2p+3<m$, one has $r<n$ and lemma \ref{t:q-rel-f} applies.
  Relations (\ref{e:q-rel-f}) can be easily translated to the ring
  $\HH_Q^*(\PPt)$ with the aid of lemma \ref{t:c-ring}, obtaining the
  two relations (\ref{e:q-rel-g}).  Then $\tilde g_1$, $\tilde g_2$
  are the classical relations (\ref{e:c-rel-g}) evaluated in the
  quantum cohomology ring.  By theorem 2.2 in \cite{vafa}, the quantum
  cohomology ring of $\HH_Q^*(\PPt)$ is the ring generated by $k$,
  $\eta$ with the two relations $\tilde g_1$, $\tilde g_2$ (compare
  also with proposition 10 of \cite{fupa}).
\end{proof}



\end{document}